\documentclass[12pt, reqno]{amsart}
\usepackage{amsmath, amstext, amsbsy, amssymb, amscd}

\newtheorem{theorem}{Theorem}[section]
\newtheorem{lemma}[theorem]{Lemma}
\newtheorem{proposition}[theorem]{Proposition}

\theoremstyle{definition}

\newtheorem{example}[theorem]{Example}

\theoremstyle{remark}
\newtheorem{remark}[theorem]{Remark}

\numberwithin{equation}{section}

\newcommand{\C}{ \mathbb C }

\newcommand{\End}{{\rm End}}

\newcommand{\Pee}{{\mathbb P}}

\newcommand{\vac}{|0\rangle}

\newcommand{\Xn}{ X^{[n]}}

\newcommand{\Z}{ \mathbb Z}

{\vskip-\lastskip\medskip
  \noindent
  {\em #1.}\enspace
  }%
{\qed\par\medskip
  }

\begin{document}
\title[Cohomology rings via Jack polynomials]
{The cohomology rings of Hilbert schemes via Jack polynomials}

\author[Wei-Ping Li]{Wei-Ping Li$^1$}
\address{Department of Mathematics, HKUST, Clear Water Bay,
  Kowloon, Hong Kong}
\email{mawpli@ust.hk}
\thanks{${}^1$Partially supported by the grant HKUST6114/02P}

\author[Zhenbo Qin]{Zhenbo Qin$^2$}
\address{Department of Mathematics, University of Missouri,
  Columbia, MO 65211, USA}
\email{zq@math.missouri.edu}
\thanks{${}^2$Partially supported by an NSF grant and
a University of Missouri Research Board grant}

\author[Weiqiang Wang]{Weiqiang Wang$^3$}
\address{Department of Mathematics, University of Virginia,
  Charlottesville, VA 22904}
\email{ww9c@virginia.edu}
\thanks{${}^3$Partially supported by an NSF grant}


\maketitle
\date{}

\section{\bf Introduction}

Fundamental and deep connections have been developed in recent
years between the geometry of Hilbert schemes $\Xn$ of points on a
(quasi-)projective surface $X$ and combinatorics of symmetric
functions. Among distinguished classes of symmetric functions, let
us mention the monomial symmetric functions, Schur polynomials,
Jack polynomials (which depend on a {\em Jack parameter}), and
Macdonald polynomials, etc (cf. \cite{Mac}). The monomial
symmetric functions can be realized as certain ordinary cohomology
classes of the Hilbert schemes associated to an embedded curve in
a surface (cf. \cite{Na1}). Nakajima \cite{Na2} further showed
that the Jack polynomials whose Jack parameter is a positive
integer $\gamma$ can be realized as certain $\mathbb
T$-equivariant cohomology classes of the Hilbert schemes of points
on the surface $X(\gamma)$ which is the total space of the line
bundle $\mathcal O_{\Pee^1}(-\gamma)$ over the complex projective
line $\Pee^1$. Here and below $\mathbb T$ stands for the
one-dimensional complex torus. In other words, the Jack parameter
is interpreted as minus the self-intersection number of the
zero-section in $X(\gamma)$. With very different motivations,
Haiman (cf. \cite{Hai} and the references therein) developed
connections between the Macdonald polynomials and the geometry of
Hilbert schemes, and in particular realized the Macdonald
polynomials as certain $\mathbb T$-equivariant $K$-homology
classes of the Hilbert schemes of points on the affine plane
$\C^2$ (A similar result has been conjectured in \cite{Na2}).

In this note, we shall establish a link somewhat different from
\cite{Na2} between equivariant cohomology of Hilbert schemes and
Jack polynomials, and then use this to describe the equivariant
and ordinary cohomology rings of the Hilbert schemes of points on
the surface $X(\gamma)$. We first show that the Jack polynomials
can be realized in terms of certain $\mathbb T$-equivariant
cohomology classes of the Hilbert schemes of points on {\em the
affine plane}, and the Jack parameter comes from the ratio of the
$\mathbb T$-weights on the two affine lines preserved by the
$\mathbb T$-action. In our view, the present construction is
conceptually simpler than the original one in \cite{Na2}. This
result is probably not very surprising however and could be well
anticipated by experts (as it is done by elaborating the ideas of
\cite{Na2} with new inputs from \cite{Vas}). But we feel it is
nice to formulate it precisely and to make it accessible to the
public. In the case when the ratio mentioned above is one, the
Jack polynomials are the Schur polynomials, and our current
construction specializes to \cite{Vas} (also cf. \cite{LQW}).

In addition, we study the $\mathbb T$-equivariant cohomology ring
of $X(\gamma)^{[n]}$ with respect to a certain $\mathbb T$-action.
The $\mathbb T$-action used in this note has the property that the
$\mathbb T$-fixed points are isolated, and is quite different from
the one used in \cite{Na2}. Generalizing an idea in \cite{Vas}, we
can define a ring structure on $H_{\mathbb
T}^{2n}(X(\gamma)^{[n]})$, which is modified from and in turn
encodes the $\mathbb T$-equivariant cohomology ring structure of
$X(\gamma)^{[n]}$. We identify the classes in $H_{\mathbb T}^{2n}
(X(\gamma)^{[n]})$ associated to the $\mathbb T$-fixed points in
$X(\gamma)^{[n]}$ with $2$-tuple Jack polynomials, by using the
connections between the Hilbert schemes of points on $\C^2$ and
Jack polynomials as formulated above. We further prove that the
ring $H_{\mathbb T}^{2n}(X(\gamma)^{[n]})$ can be simply described
by declaring these $2$-tuple Jack polynomials to be (essentially)
idempotents.

In both cases when $X$ is $\C^2$ or $X(\gamma)$, our main
technical tool is the construction of a Heisenberg algebra acting
on the direct sum $\displaystyle{\bigoplus_{n=0}^\infty H_{\mathbb
T}^{2n}(\Xn)}$ of the equivariant cohomology groups of middle
degree. This construction generalizes the one in \cite{Vas}, which
is in turn a modification of the construction in \cite{Na1} (also
cf. \cite{Gro}).

Finally, we note that the ordinary cohomology ring of the Hilbert
scheme $X(\gamma)^{[n]}$ can be shown to be isomorphic to the
graded ring associated to a natural filtration on the ring
$H_{\mathbb T}^{2n}(X(\gamma)^{[n]})$. In this way, we obtain an
algorithm for computing the ordinary cup product of cohomology
classes in $X(\gamma)^{[n]}$.

\section{\bf Equivariant cohomology rings of Hilbert schemes}
\subsection{Surfaces with torus actions}
\par
$\,$

Let ${\mathbb T}= \mathbb C^*$, and $\theta$ be the
$1$-dimensional standard ${\mathbb T}$-module. For an algebraic
variety $M$ with a ${\mathbb T}$-action, let $H^*_{\mathbb T}(M)$
be the equivariant cohomology ring of $M$ with $\mathbb C$
coefficients. Then $H^*_{\mathbb T}(M)$ is a $\mathbb C[t]$-module
if we identify $H^*_{\mathbb T}(\text{\rm pt})$ and $\mathbb C[t]$
($t$ is an element of degree-$2$). For a ${\mathbb T}$-equivariant
and proper morphism $f\colon N\to M$ of algebraic varieties, there
is a Gysin homomorphism $f_!\colon H^*_{\mathbb T}(N)\to
H_{\mathbb T}^*(M)$ of ${\mathbb T}$-equivariant cohomology
groups. If $N$ is a ${\mathbb T}$-equivariant codimension-$k$
closed subvariety of $M$ and $i\colon N\to M$ is the inclusion
map, define $[N]=i_!(1_N) \in H^{2k}_{\mathbb T}(M)$ where $1_N
\in H^0_{\mathbb T}(N)$ is the unit of the algebra $H^*_{\mathbb
T}(N)$.

In this note, we shall consider two types of surfaces with
${\mathbb T}$-actions, and the Hilbert schemes of points on these
surfaces. The first one is the complex plane $\mathbb C^2$, while
the second is the total space of a line bundle over $\mathbb P^1$.
The ${\mathbb T}$-actions on these surfaces are specified in the
following two examples.

\begin{example} \label{complexplane}
Fix two nonzero integers $\alpha$ and $\beta$ with the same signs.
Let $u, v$ be the standard coordinate functions on $\mathbb C^2$.
We define the action of $\mathbb T$ on $\mathbb C^2$ by
\begin{eqnarray}  \label{action_plane}
s\cdot (u, v)=(s^{\alpha}u, s^{-\beta}v), \quad s\in \mathbb T.
\end{eqnarray}
The origin of $\mathbb C^2$ is the only fixed point, which will be
denoted by $x$. Let $\Sigma$ and $\Sigma^{\prime}$ be the $u$-axis
and $v$-axis respectively in $\mathbb C^2$. As $\mathbb
T$-modules, we have $T_{x}\Sigma = \theta^{-\alpha}$ and
$T_{x}\Sigma^{\prime}=\theta^{\beta}$. By the localization theorem
(see \cite{C-K}), we get
\begin{eqnarray}  \label{Sigma_x}
[\Sigma]=-{\alpha}^{-1} t^{-1}[x], \quad
[\Sigma^{\prime}]={\beta}^{-1} t^{-1}[x].
\end{eqnarray}
\end{example}

\begin{example} \label{linebundle}
Fix an integer $\gamma>1$. Let $X(\gamma)$ be the total space of
the line bundle $\mathcal O_{\mathbb P^1}(-\gamma)$ over $\Pee^1$.
The quasi-projective surface $X(\gamma)$ can be regarded as the
quotient space of $\mathbb C \times (\mathbb C^2-\{0\})$ by the
$\mathbb C^*$-action defined by
\begin{eqnarray}   \label{equi_cls}
s\cdot (b,b_1, b_2)=( s^{-\gamma}b, sb_1, sb_2), \quad s \in
\mathbb C^*.
\end{eqnarray}
We use $[(b,b_1,b_2)]$ to denote the equivalence class. Define a
$\mathbb T$-action on $X(\gamma)$ by
\begin{eqnarray} \label{action_X}
s\cdot [(b,b_1, b_2)]=[(sb,s^{-1}b_1, b_2)], \quad s \in\mathbb T.
\end{eqnarray}
For $i=1, 2$, let $X_i$ be the open subset of $X(\gamma)$ given by
\begin{eqnarray} \label{open-cover}
X_1=\{[(b,b_1, b_2)]\,|\, b_{2}=1\},\quad X_2=\{[(b,b_1,
b_2)]\,|\, b_{1}=1\}.
\end{eqnarray}
Then $X_1$ and $X_2$ form an affine open cover of $X(\gamma)$.
Moreover, each $X_i$ is $\mathbb T$-invariant. For simplicity,
denote the point $[(b,b_1,1)] \in X_1$ by $(b,b_1)$. Similarly,
denote $[(b,1,b_2)] \in X_2$ by $(b,b_2)$. Then $\mathbb T$ acts
on the points of $X_1$ by $s\cdot (b,b_1)=(sb,s^{-1}b_1)$, i.e.,
$\mathbb T$ acts on the coordinate functions $u_1$ and $v_1$ of
$X_1$ by
\begin{eqnarray} \label{action_X1}
s\cdot (u_1, v_1)=(s^{-1}u_1, sv_1), \quad s \in\mathbb T.
\end{eqnarray}
Similarly, $\mathbb T$ acts on the coordinate functions $u_2$ and
$v_2$ of $X_2$ by
\begin{eqnarray} \label{action_X2}
s\cdot (u_2, v_2)=(s^{\gamma-1}u_2, s^{-1}v_2), \quad s \in\mathbb
T.
\end{eqnarray}

Let $x_i$ be the origin of $X_i$. Then $X(\gamma)^{\mathbb T}
=\{x_1,x_2\}$. Let $\rho: X(\gamma) \to \mathbb P^1$ be the
projection sending $[(b,b_1, b_2)]$ to $[b_1, b_2]$. Let $\Sigma_0
\cong \mathbb P^1$ be the zero section of $\rho$, and
\begin{eqnarray} \label{sigma12}
\Sigma_1 = \rho^{-1}([0,1]), \qquad \Sigma_2=\rho^{-1}([1,0]).
\end{eqnarray}
Then as $\mathbb T$-modules, we have $T_{x_1}\Sigma_1=\theta, \,\,
T_{x_1}\Sigma_0 = \theta^{-1},\,\, T_{x_2}\Sigma_0 =\theta$, and
$T_{x_2}\Sigma_2 = \theta^{1-\gamma}$. By the localization
theorem, we get
\begin{eqnarray} \label{Sigma12}
[\Sigma_1]=t^{-1}[x_1],\,\,\, [\Sigma_0]=-t^{-1}[x_1]+
t^{-1}[x_2], \,\,\, [\Sigma_2]= (1-\gamma)^{-1} t^{-1}[x_2].
\end{eqnarray}
\end{example}

\begin{remark} (i) Let $X(1)$ be the total space of
$\mathcal O_{\mathbb P^1}(-1)$. Using $[(b,b_1,b_2)] \in X(1)$ to
denote the equivalence class defined by (\ref{equi_cls}), we
define a $\mathbb T$-action on $X(1)$ by
\begin{eqnarray*}
s \cdot [(b,b_1, b_2)]=[(sb,s^{-2}b_1, b_2)], \quad s \in \mathbb
T.
\end{eqnarray*}
Then our methods and results below apply to $X(1)$ as well.

(ii) Other $\mathbb T$-actions on the surfaces $X(\gamma)$,
$\gamma \ge 1$, with isolated fixed points can be treated
similarly.
\end{remark}
\subsection{Distinguished equivariant cohomology classes}
\par
$\,$

In the rest of this note, let $X$ be a surface in
Example~\ref{complexplane} or Example~\ref{linebundle}. Our goal
in this subsection is to define some distinguished equivariant
cohomology classes for the Hilbert schemes. Let $\Xn$ be the
Hilbert scheme parametrizing all the $0$-dimensional closed
subschemes $\xi$ of $X$ with $\dim_{\mathbb C} H^0(\mathcal
O_{\xi})=n$. The $\mathbb T$-action on $X$ induces a $\mathbb
T$-action on $\Xn$. The support of a ${\mathbb T}$-fixed point in
$\Xn$ is contained in $X^{\mathbb T}$. By the results in
\cite{E-S}, the $\mathbb T$-fixed points of $\Xn$ are isolated and
parametrized in terms of (multi-)partitions.


Next, let $X = \C^2$ as in Example~\ref{complexplane}. The
$\mathbb T$-fixed points of $\Xn$ are supported in $X^T=\{x\}$ and
indexed by partitions $\lambda$ of $n$. We use $\xi_{\lambda}$ to
denote the fixed point in $(X^{[n]})^{\mathbb T}$ corresponding to
a partition $\lambda$ of $n$. Then for $\lambda \vdash n$, the
tangent space of $\Xn$ at the fixed point $\xi_{\lambda}$ is
$\mathbb T$-equivariantly isomorphic to (see \cite{E-S, Na1,
Na2}):
\begin{eqnarray}  \label{tangent2}
T_{\xi_{\lambda}}\Xn = \bigoplus_{\square \in D_{\lambda}} \left (
\theta^{\alpha(\ell(\square)+1)+\beta a(\square)} \oplus
\theta^{-\alpha\ell(\square)- \beta (a(\square)+1)} \right )
\end{eqnarray}
where $D_{\lambda}$ is the Young diagram associated to the
partition $\lambda$, $\square$ is a cell in $D_{\lambda}$,
$\ell(\square)$ is the leg of $\square$, and $a(\square)$ is the
arm of $\square$ (see \cite{Mac} for the notations). So
\begin{eqnarray}  \label{eTC2}
e_{\mathbb T}(T_{\xi_{\lambda}}\Xn)=(-1)^n c_{\lambda}(\alpha,
\beta) c'_{\lambda}(\alpha, \beta) t^{2n}.
\end{eqnarray}
where $e_{\mathbb T}( \cdot )$ stands for the equivariant Euler
class and
\begin{eqnarray}
c_{\lambda}(\alpha,\beta)&=&\prod_{\square\in D_{\lambda}}
  \big(\alpha(\ell(\square)+1)+\beta a(\square)\big),
  \label{def_c} \\
c'_{\lambda}(\alpha,\beta)&=&\prod_{\square\in
  D_{\lambda}}\big(\alpha\ell(\square)+ \beta (a(\square)+1)
  \big).  \label{def_c'}
\end{eqnarray}
Note that $[\xi_{\lambda}] \in H^{4n}_{\mathbb T}(\Xn)$. We define
the following distinguished class:
\begin{eqnarray} \label{def_[la]}
[\lambda]=\frac{(-1)^n}{c_{\lambda}(\alpha, \beta)}t^{-n}
[\xi_{\lambda}].
\end{eqnarray}

Now let $X = X(\gamma)$ as in Example~\ref{linebundle}. In view of
(\ref{action_plane}) and (\ref{action_X1}), there is a $\mathbb
T$-equivariant identification between $X_1$ and the complex plane
in Example~\ref{complexplane} with $\alpha = \beta = -1$;
similarly for $X_2$ and the complex plane in
Example~\ref{complexplane} with $\alpha = \gamma-1$ and $\beta=
1$. Then the $\mathbb T$-fixed points of $\Xn$ are of the form
$\xi_{\lambda^1} + \xi_{\lambda^2}$ where $\lambda^1$ and
$\lambda^2$ are partitions with $|\lambda^1| + |\lambda^2|= n$,
and $\xi_{\lambda^1}$ and $\xi_{\lambda^2}$ are defined in the
previous paragraph as we identify $X_1$ and $X_2$ with $\C^2$. For
simplicity, put
\begin{eqnarray*}
\xi_{\lambda^1, \lambda^2} = \xi_{\lambda^1} + \xi_{\lambda^2}.
\end{eqnarray*}
We have a $\mathbb T$-equivariant splitting $T_{\xi_{\lambda^1,
\lambda^2}} \Xn \cong T_{\xi_{\lambda^1}} X_1^{[|\lambda^1|]}
\oplus T_{\xi_{\lambda^2}} X_2^{[|\lambda^2|]}$. By (\ref{eTC2}),
\begin{eqnarray*}  \label{eTlinebdle}
e_{\mathbb T}(T_{\xi_{\lambda^1, \lambda^2}} \Xn) =(-1)^n
c_{\lambda^1}(-1, -1) c'_{\lambda^1}(-1, -1)
c_{\lambda^2}(\gamma-1, 1) c'_{\lambda^2}(\gamma-1, 1) t^{2n}.
\end{eqnarray*}
Also, as in (\ref{def_[la]}), we introduce the distinguished class
\begin{eqnarray} \label{def_[la12]}
[\lambda^1, \lambda^2]= \frac{(-1)^n}{c_{\lambda^1}(-1, -1)
c_{\lambda^2}(\gamma-1, 1)}t^{-n}[\xi_{\lambda^1, \lambda^2}].
\end{eqnarray}

\subsection{Bilinear forms on the equivariant cohomology}
\label{subsect_form}
\par
$\,$

It is known from \cite{E-S, Got} that the odd Betti numbers of
$\Xn$ are equal to zero and $H^k(\Xn)=0$ for $k>2n$. Hence the
spectral sequence associated with the fibration
$\Xn\times_{\mathbb T}E{\mathbb T}\to B{\mathbb T}$ degenerates at
the $E_2$-term. We have
\begin{eqnarray*}
H^{2k}_{\mathbb T}(\Xn) =t^{k-n}\cup H_{\mathbb T}^{2n}(\Xn)
\end{eqnarray*}
for $k\ge n$. Therefore, the classes defined in (\ref{def_[la]})
and (\ref{def_[la12]}) are contained in
\begin{eqnarray*}
\mathbb H_n \stackrel{\text{def}}{=} H^{2n}_{\mathbb T}(\Xn).
\end{eqnarray*}
Moreover, we can define a product structure $\star$ on $\mathbb
H_n$ as follows (also cf. \cite{Vas}):
\begin{eqnarray}   \label{star}
t^n\cup(A\star B)=A\cup B \in H^{4n}_{\mathbb T}(\Xn)
\end{eqnarray}
for $A,\, B\in H^{2n}_{\mathbb T}(\Xn)$. We see that
$(H^{2n}_{\mathbb T}(\Xn),\star)$ is a ring.

Let $H^*_{\mathbb T}(\cdot)^{\prime} = H^*_{\mathbb
T}(\cdot)\otimes_{\mathbb C[t]} \mathbb C(t)$ be the localization,
and let
\begin{eqnarray*}
\iota: (\Xn)^{\mathbb T}\to \Xn
\end{eqnarray*}
be the inclusion map. By abusing notations, we also use $\iota_!$
to denote the induced Gysin map on the localized equivariant
cohomology groups:
\begin{eqnarray}\label{Gysin}
\iota_!\colon H^*_{\mathbb T}((\Xn)^{\mathbb T})^{\prime}\to
H^*_{\mathbb T}(\Xn)^{\prime},
\end{eqnarray}
which is an isomorphism by the localization theorem.

Define a bilinear form $\langle -, -\rangle: H^*_{\mathbb
T}(\Xn)^{\prime} \times H^*_{\mathbb T}(\Xn)^{\prime} \to \mathbb
C(t)$:
\begin{eqnarray} \label{bilinear1}
\langle A, B \rangle =(-1)^np_!\iota_!^{-1}(A\cup B)
\end{eqnarray}
where $p$ is the projection $(\Xn)^{\mathbb T}\to \text{pt}$. This
induces a bilinear form $\langle -, -\rangle$ on
\begin{eqnarray} \label{def_HX'}
\mathbb H_X' \,\,\, \stackrel{\text{\rm def}}{=} \,\,\,
\bigoplus_{n=0}^{\infty} H^*_{\mathbb T}(\Xn)^{\prime}.
\end{eqnarray}

Next, we study the restriction of the bilinear form $\langle -,
-\rangle$ to $\mathbb H_n = H^{2n}_{\mathbb T}(\Xn)$. When $X =
\C^2$, we see from the projection formula and (\ref{eTC2}) that
\begin{eqnarray} \label{cup-prod}
   [\xi_{\lambda}] \cup[\xi_{\mu}]
&=&i_{\lambda!}(1_{\xi_{\lambda}})
   \cup i_{\mu!}(1_{\xi_{\mu}}) =i_{\lambda!}\big (1_{\xi_{\lambda}}
   \cup i^*_{\lambda}i_{\mu!}(1_{\xi_{\mu}}) \big )  \nonumber \\
&=&\delta_{\lambda,\mu}e_{\mathbb T}(T_{\xi_{\lambda}}\Xn)
   [\xi_{\lambda}] = \delta_{\lambda,\mu}(-1)^n
   c_{\lambda}(\alpha, \beta) c'_{\lambda}(\alpha, \beta)
   t^{2n} [\xi_{\lambda}]. \qquad
\end{eqnarray}
It follows from (\ref{def_[la]}) and (\ref{bilinear1}) that for
$\lambda, \mu \vdash n$, we have
\begin{eqnarray} \label{pairing_la_mu1}
\langle[\lambda],[\mu]\rangle =\delta_{\lambda,\mu}
\frac{c'_{\lambda}(\alpha,\beta)}{c_{\lambda}(\alpha,\beta)}.
\end{eqnarray}
By the localization theorem, we see that the classes $[\lambda]$,
$\lambda \vdash n$ form a linear basis of the $\mathbb C$-vector
space $\mathbb H_n$. Similarly, when $X = X(\gamma)$ is from
Example~\ref{linebundle},
\begin{eqnarray} \label{pairing_la_mu2}
\langle[\lambda^1, \lambda^2], [\mu^1,\mu^2]\rangle =
\delta_{\lambda^1,\mu^1}\delta_{\lambda^2,\mu^2}
\frac{c'_{\lambda^2}(1-\gamma, -1)}{c_{\lambda^2}(1-\gamma,-1)} =
\langle[\lambda^1],[\mu^1]\rangle \cdot
\langle[\lambda^2],[\mu^2]\rangle,
\end{eqnarray}
and the classes $[\lambda^1, \lambda^2]$, where $|\lambda^1| +
|\lambda^2| = n$, form a linear basis of $\mathbb H_n$.

It follows that the restriction to $\mathbb H_n$ of the bilinear
form $\langle -, -\rangle$ on $H^*_{\mathbb T}(\Xn)^{\prime}$ is a
nondegenerate bilinear form $\langle -, -\rangle: \mathbb H_n
\times \mathbb H_n \to \C$. This induces a nondegenerate bilinear
form $\langle -, -\rangle: \mathbb H_X \times \mathbb H_X \to \C$
where the space $\mathbb H_X$ is defined by
\begin{eqnarray} \label{def_HX}
\mathbb H_X \,\,\, = \,\,\, \bigoplus_{n=0}^{\infty} \mathbb H_n.
\end{eqnarray}

\section{\bf Heisenberg algebras, equivariant cohomology and
Jack polynomials}
\subsection{Heisenberg algebra actions}
\label{general}
\par
$\,$

Let $X$ be a surface in Example~\ref{complexplane} or
Example~\ref{linebundle}. Fix a positive integer $i$. For a
$\mathbb T$-invariant closed curve $Y\subset X$, we define
\begin{eqnarray*}
Y_{n, i}=\{(\xi, \eta)\in X^{[n+i]}\times \Xn\, |\, \eta\subset
\xi, \,\, \text{Supp}(I_{\eta}/I_{\xi})=\{y\}, \,\, y \in Y\}
\end{eqnarray*}
where $I_{\eta}$ and $I_{\xi}$ are the sheaves of ideals
corresponding to $\eta$ and $\xi$ respectively. Let $p_1$ and
$p_2$ be the projections of $ X^{[n+i]}\times \Xn$ to the two
factors. As in \cite{Vas}, we define a linear operator $\mathfrak
p_{-i}([Y])\in \End(\mathbb H_X^{\prime})$ by
\begin{eqnarray}\label{creation}
\mathfrak p_{-i}([Y])(A)=p_{1!}(p_2^*A\cup [Y_{n, i}])
\end{eqnarray}
for $A\in H^*_{\mathbb T}(\Xn)^{\prime}$. Note that the
restriction of $p_1$ to $Y_{n, i}$ is proper. We define $\mathfrak
p_i([Y])\in \End(\mathbb H_X^{\prime})$ to be the adjoint operator
of $\mathfrak p_{-i}([Y])$ with respect to the bilinear form
$\langle-, -\rangle$ on $\mathbb H_X^{\prime}$. For $A\in
H^*_{\mathbb T}(\Xn)^{\prime}$, we have
\begin{eqnarray}  \label{deletion}
\mathfrak p_{i}([Y])(A)=(-1)^ip^{\prime}_{2!}(\iota\times
\text{Id})_!^{-1}(p_1^*A\cup [Y_{n-i, i}])
\end{eqnarray}
where $p_2^{\prime}$ is the projection of $(\Xn)^{\mathbb T}\times
X^{[n-i]}$ to $X^{[n-i]}$. Finally we put $\mathfrak p_0([Y])=0$.

By the definition of $\mathfrak p_{-i}([Y])$ for $i>0$, its
restriction to $\mathbb H_X$ gives a linear operator in
$\End(\mathbb H_X)$, denoted by $\mathfrak p_{-i}([Y])$ as well.
Next, we recall from Subsection~\ref{subsect_form} that there is a
nondegenerate bilinear form $\langle -, -\rangle\colon \mathbb
H_X\otimes_{\mathbb C} \mathbb H_X \to \mathbb C$, which is the
restriction of the bilinear form $\langle -, -\rangle$ on $\mathbb
H_X^{\prime}$. Thus, the restriction of $\mathfrak p_i([Y])$ to
$\mathbb H_X$ is the adjoint operator of $\mathfrak p_{-i}([Y])$
with respect to the bilinear form $\langle-, -\rangle$ on $\mathbb
H_X$, and hence is an operator in $\End(\mathbb H_X)$ which will
again be denoted by $\mathfrak p_i([Y])$.

By (\ref{Sigma_x}) and (\ref{Sigma12}), $H^2_{\mathbb T}(X)$ is
linearly spanned by the classes $[Y]$ where $Y$ denotes $\mathbb
T$-invariant closed curves in $X$. So we can extend the notion
$\mathfrak p_k([Y])$ linearly to obtain the operator $\mathfrak
p_k(\omega) \in \End(\mathbb H_X)$ for an arbitrary class $\omega
\in H^2_{\mathbb T}(X)$. Note from Example~\ref{complexplane}
(respectively, Example~\ref{linebundle}) that $\Sigma$ and
$\Sigma'$ (respectively, $\Sigma_0$ and $\Sigma_i$ where $i = 1$
or $2$) intersect transversely at one point. This (simple but
crucial) observation together with an argument parallel to
\cite{Na1, Vas} leads to the following.

\begin{proposition} \label{prop1}
The operators $\mathfrak p_k(\omega)$, $k\in \mathbb Z$ and
$\omega \in \mathbb H_1 = H^2_{\mathbb T}(X)$, acting on $\mathbb
H_X$ satisfy the following Heisenberg commutation relation:
\begin{eqnarray}\label{Heisenberg}
[\mathfrak p_k(\omega_1), \mathfrak p_{\ell}(\omega_2)]=
k\delta_{k, -\ell}\langle \omega_1, \omega_2\rangle \,\, \text{\rm
Id}.
\end{eqnarray}
Furthermore, $\mathbb H_X$ becomes the Fock space over the
Heisenberg algebra modeled on $H^2_{\mathbb T}(X)$ with the unit
$|0\rangle \in H^0_{\mathbb T}(X^{[0]})$ of $H^*_{\mathbb
T}(X^{[0]})$ being a highest weight vector. \qed
\end{proposition}
\subsection{The case of the complex plane $\mathbb C^2$}
\label{sub_plane}
\par
$\,$

In this subsection, we consider $X=\mathbb C^2$ from
Example~\ref{complexplane}. Let $\mathfrak p_i = \mathfrak
p_i([\Sigma])$ for $i\in \mathbb Z$. Given a partition
$\lambda=(\lambda_1, \ldots, \lambda_\ell) =
(1^{m_1}2^{m_2}\ldots)$, define
\begin{eqnarray} \label{frakp_la}
\mathfrak z_{\lambda}&=&\prod_{i\ge 1}r^{m_i}m_i!, \nonumber \\
\mathfrak p_{-\lambda}&=&\frac{1}{\mathfrak z_{\lambda}} \,\,
\prod_{i\ge 1}\mathfrak p_{-i}^{m_i}.
\end{eqnarray}
By (\ref{Sigma_x}), (\ref{def_[la]}) and (\ref{pairing_la_mu1}),
we have $\langle[\Sigma], [\Sigma]\rangle={\beta}/{\alpha}$. So by
Proposition~\ref{prop1},
\begin{eqnarray} \label{pairing_frakp_la_mu}
\langle\mathfrak p_{-\lambda}|0\rangle, \mathfrak
p_{-\mu}|0\rangle\rangle = \delta_{\lambda,\mu} \,\,
\frac{1}{\mathfrak z_{\lambda}} \,\,
({\beta}/{\alpha})^{\ell(\lambda)}.
\end{eqnarray}

Let $S^nX$ (respectively, $S^n\Sigma$) be the $n$-th symmetric
product of $X$ (respectively, of $\Sigma$). Let $\pi\colon \Xn \to
S^nX$ be the Hilbert-Chow morphism. For $\lambda \vdash n$, define
\begin{eqnarray*}
S^n_{\lambda}\Sigma = \left \{\sum_{i=1}^\ell \lambda_i y_i \in
S^n\Sigma\,|\, y_i \in \Sigma, \,\, \text{\rm and } y_i \neq y_j
\,\, \text{if } i\neq j \right \}.
\end{eqnarray*}
Let $\Sigma^{(\lambda)}$ be the closure of
$\pi^{-1}(S^n_{\lambda}\Sigma)$ in $\Xn$. An argument parallel to
the proof of the Corollary 6.10 in \cite{Na2} (also cf.
\cite{Vas}) shows that
\begin{eqnarray} \label{sym_[la]}
[\Sigma^{(\lambda)}] = [\lambda]+\sum_{\mu<\lambda} c_{\lambda,
\mu}[\mu]
\end{eqnarray}
where $c_{\lambda, \mu} \in \C$ and ``$<$" denotes the dominance
partial ordering of partitions.

Let $\Lambda$ be the ring of symmetric polynomials in infinitely
many variables, $\Lambda^n$ be the space of degree-$n$ symmetric
polynomials, $m_{\lambda}$ be the monomial symmetric function
associated to a partition $\lambda$, and $p_k$ be the $k$-th power
sum symmetric function. Given a partition
$\lambda=(1^{m_1}2^{m_2}\ldots)$, we define
\begin{eqnarray} \label{p_la}
p_{\lambda}= \displaystyle\frac{1}{\mathfrak z_{\lambda}}
\prod_{i\ge 1}p_{i}^{m_i}.
\end{eqnarray}
It is known that the symmetric functions $p_{\lambda}$ form a
linear basis of the $\C$-vector space $\Lambda$. So we can define
a bilinear form $\langle-, -\rangle$ on $\Lambda$ by
\begin{eqnarray} \label{pairing_p_la_mu}
\langle p_{\lambda}, p_{\mu}\rangle = \delta_{\lambda,\mu} \,\,
\frac{1}{\mathfrak z_{\lambda}} \,\,
({\beta}/{\alpha})^{\ell(\lambda)}.
\end{eqnarray}
Let $P_{\lambda}^{({\beta}/{\alpha})}$ denote the Jack polynomials
(see page 379 in \cite{Mac} for their characterization). We
introduce a ring structure $\circ$ on $\Lambda^n$ defined by
\begin{eqnarray}  \label{prod}
\frac{P_{\lambda}^{({\beta}/{\alpha})}}
  {c^{\prime}_{\lambda}(\alpha,\beta)}
\circ \frac{P_{\mu}^{({\beta}/{\alpha})}}
  {c^{\prime}_{\mu}(\alpha,\beta)}
=\delta_{\lambda,\mu}\frac{P_{\lambda}^{({\beta}/{\alpha})}}
  {c^{\prime}_{\lambda}(\alpha,\beta)}.
\end{eqnarray}

\begin{theorem} \label{thm1}
There exists a linear isomorphism $\phi: \mathbb H_X \to \Lambda$
preserving bilinear forms such that $\phi(\mathfrak
p_{-\lambda}|0\rangle)=p_{\lambda}$,
$\phi([\Sigma^{(\lambda)}])=m_{\lambda}$, and $\phi([\lambda])
=P_{\lambda}^{({\beta}/{\alpha})}$. Furthermore, the restriction
$\phi_n$ of $\phi$ to $\mathbb H_n$ is an isomorphism of rings,
i.e., $\phi_n: (\mathbb H_n, \,\star) \cong (\Lambda^n, \circ)$.
\end{theorem}
\begin{proof}
The linear isomorphism $\phi: \mathbb H_X \to \Lambda$ is defined
by mapping $\mathfrak p_{-\lambda}|0\rangle$ to $p_{\lambda}$. So
\begin{eqnarray} \label{phi_commu_p}
\phi(\mathfrak p_{-i}A)= p_i \,\, \phi(A), \qquad i > 0, \,\,\, A
\in \mathbb H_X.
\end{eqnarray}

Next, an argument similar to the proof of the Theorem 4.6 in
\cite{Na2} (also cf. \cite{Vas}) verifies that for all $i > 0$, we
have
\begin{eqnarray} \label{ps*m}
\mathfrak p_{-i}\cdot [\Sigma^{(\lambda)}]=\sum_\mu a_{\lambda,
\mu}[\Sigma^{(\mu)}]
\end{eqnarray}
where the coefficients $a_{\lambda, \mu}$ are the same as in $p_i
\,\, m_{\lambda}=\sum_{\mu}a_{\lambda, \mu}m_{\mu}$.
Therefore we conclude from an induction, (\ref{ps*m}) and
(\ref{phi_commu_p}) that $\phi([\Sigma^{(\lambda)}])=
m_{\lambda}$.

By (\ref{pairing_frakp_la_mu}) and (\ref{pairing_p_la_mu}), $\phi$
preserves the bilinear forms. So we see from
(\ref{pairing_la_mu1}) that
\begin{eqnarray*}
  \langle \phi([\lambda]), \phi([\mu]) \rangle
= \delta_{\lambda,\mu} \frac{c'_{\lambda}(\alpha,\beta)}
  {c_{\lambda}(\alpha,\beta)}
= \delta_{\lambda,\mu} \frac{c'_{\lambda}(1,\beta/\alpha)}
  {c_{\lambda}(1,\beta/\alpha)}.
\end{eqnarray*}
By (\ref{sym_[la]}), we have $[\lambda]=[\Sigma^{(\lambda)}]
+\sum_{\mu<\lambda} d_{\lambda, \mu}[\Sigma^{(\mu)}]$. It follows
that
\begin{eqnarray*}
\phi([\lambda]) = m_{\lambda} + \sum_{\mu<\lambda} d_{\lambda,
\mu} m_{\mu}.
\end{eqnarray*}
By the characterization of the Jack polynomials in \cite{Mac},
$\phi([\lambda]) = P_{\lambda}^{({\beta}/{\alpha})}$.

Finally, by (\ref{cup-prod}) and the definition of $\star$, we get
\begin{eqnarray*}
\displaystyle{\frac{[\lambda]}
{c^{\prime}_{\lambda}(\alpha,\beta)} \star
\frac{[\mu]}{c^{\prime}_{\mu}(\alpha,\beta)}
=\delta_{\lambda,\mu}\frac{[\lambda]}
{c^{\prime}_{\lambda}(\alpha,\beta)}}, \qquad \lambda, \mu \vdash
n.
\end{eqnarray*}
In view of (\ref{prod}), $\phi_n: (\mathbb H_n, \,\star) \to
(\Lambda^n, \circ)$ is an isomorphism of rings.
\end{proof}
\subsection{The case of the surface $X(\gamma)$}
\label{sub:line-bundle}
\par
$\,$

Let $X = X(\gamma)$ be the surface from Example~\ref{linebundle}.
Recall that there is a $\mathbb T$-equivariant identification
between the affine open subset $X_1$ of $X$ and the complex plane
in Example~\ref{complexplane} with $\alpha = \beta = -1$;
similarly for $X_2$ and the complex plane in
Example~\ref{complexplane} with $\alpha = \gamma-1$ and $\beta=
1$. By the discussions in Subsection~\ref{general}, $\mathbb
H_{X}$ is the Fock space of the Heisenberg algebra generated by
$\mathfrak p_i([\Sigma_1])$ and $\mathfrak p_i([\Sigma_2])$ with
$i \in \Z$, where $\Sigma_1$ and $\Sigma_2$ are the two $\mathbb
T$-equivariant fibers in ${X}$ defined by (\ref{sigma12}). Note
also that $\mathbb H_{X_1}$ is the Fock space of the Heisenberg
algebra generated by $\mathfrak p_i([\Sigma_1])$ for $i \in \Z$,
where $\Sigma_1$ is considered as a $\mathbb T$-equivariant closed
curve in the affine open subset $X_1 \subset {X}$. Similarly,
$\mathbb H_{X_2}$ is the Fock space of the Heisenberg algebra
generated by $\mathfrak p_i([\Sigma_2])$ for $i \in \Z$. To avoid
confusions, we use $\mathfrak p_i^{X_1}$ to denote the operators
$\mathfrak p_i([\Sigma_1])$ acting on $\mathbb H_{X_1}$,
$\mathfrak p_i^{X_2}$ to denote the operators $\mathfrak
p_i([\Sigma_2])$ acting on $\mathbb H_{X_2}$, and $\mathfrak
p_i([\Sigma_1])$ and $\mathfrak p_i([\Sigma_2])$ for the
Heisenberg operators acting on $\mathbb H_{X}$. As in
(\ref{frakp_la}), for a partition $\lambda =
(1^{m_1}2^{m_2}\ldots)$ and for $j = 1$ or $2$, we define
\begin{eqnarray*}
\mathfrak p_{-\lambda}([\Sigma_j]) = \frac{1}{\mathfrak
z_{\lambda}} \,\, \prod_{i\ge 1}\mathfrak
p_{-i}([\Sigma_j])^{m_i}.
\end{eqnarray*}

Define a linear map $\Psi: \mathbb H_{X_1} \otimes_{\mathbb C}
\mathbb H_{X_2} \to \mathbb H_X$ as follows:
\begin{eqnarray}  \label{Psi}
\Psi([\lambda^1] \otimes [\lambda^2])=[\lambda^1,\lambda^2].
\end{eqnarray}
The linear map $\Psi$ is an isomorphism since the classes
$[\lambda^1], [\lambda^2], [\lambda^1,\lambda^2]$ defined in
(\ref{def_[la]}) and (\ref{def_[la12]}) form linear bases of
$\mathbb H_{X_1}, \mathbb H_{X_2}, \mathbb H_{X}$ respectively. In
addition, we see from (\ref{pairing_la_mu2}) that $\Psi$ preserves
the bilinear forms.

\begin{lemma}\label{lem2}
The linear isomorphism $\Psi$ commutes with Heisenberg operators,
i.e.,
\begin{eqnarray*}
  \Psi \circ \big (\mathfrak p_{-i}^{X_1}\otimes \text{\rm Id}
  \big )
&=&\mathfrak p_{-i}([\Sigma_1]) \circ \Psi,       \\
  \Psi \circ \big (\text{\rm Id} \otimes \mathfrak p_{-i}^{X_2}
  \big )
&=&\mathfrak p_{-i}([\Sigma_2]) \circ \Psi.
\end{eqnarray*}
\end{lemma}
\begin{proof}
By the symmetry between $\Sigma_1$ and $\Sigma_2$, we need only to
prove the first identity. Also, since $\Psi$ preserves the
bilinear forms and $\mathfrak p_i$ is the adjoint operator of
$\mathfrak p_{-i}$, it suffices to prove the first identity for
$i>0$.

Since the two fibers $\Sigma_1$ and $\Sigma_2$ do not intersect,
we see that for partitions $\lambda^1$ and $\lambda^2$ with
$|\lambda^1| + |\lambda^2| = n$, the $\mathbb T$-equivariant
closed subvariety $\Sigma_1^{(\lambda^1)}\times
\Sigma_2^{(\lambda^2)} \subset \Xn$ is the closure of $\pi^{-1}
\big ( S^{|\lambda^1|}_{\lambda^1}\Sigma_1 \times
S^{|\lambda^2|}_{\lambda^2}\Sigma_2 \big )$ in $\Xn$, where
$\pi\colon \Xn \to S^nX$ denotes the Hilbert-Chow morphism. We
conclude that
\begin{eqnarray*}
\Psi \big ([\Sigma_1^{(\lambda^1)}]\otimes
[\Sigma_2^{(\lambda^2)}] \big ) =[\Sigma_1^{(\lambda^1)}\times
\Sigma_2^{(\lambda^2)}]
\end{eqnarray*}
by writing $[\Sigma_1^{(\lambda^1)} \times
\Sigma_2^{(\lambda^2)}]$ in terms of $[\lambda^1, \lambda^2]$ and
$[\mu^1, \mu^2]$ where $\mu^1 < \lambda^1$ or $\mu^2 < \lambda^2$,
similarly as in (\ref{sym_[la]}). It implies that the classes
$[\Sigma_1^{(\lambda^1)} \times \Sigma_2^{(\lambda^2)}]$ form a
linear basis of $\mathbb H_X$. By (\ref{ps*m}) and a similar
computation for $\mathfrak
p_{-i}([\Sigma_1])[\Sigma_1^{(\lambda^1)} \times
\Sigma_2^{(\lambda^2)}]$, we have
\begin{eqnarray*}
   \Psi \left ( \big (\mathfrak p_{-i}^{X_1}\otimes
   \text{\rm Id} \big ) \big ([\Sigma_1^{(\lambda^1)}]
   \otimes [\Sigma_2^{(\lambda^2)}] \big ) \right )
&=&\Psi \left ( \big (\mathfrak p_{-i}^{X_1}
   [\Sigma_1^{(\lambda^1)}] \big ) \otimes
   [\Sigma_2^{(\lambda^2)}]  \right )          \\
&=&\mathfrak p_{-i}([\Sigma_1])[\Sigma_1^{(\lambda^1)}
   \times \Sigma_2^{(\lambda^2)}]     \\
&=&\mathfrak p_{-i}([\Sigma_1]) \,\, \Psi \big (
   [\Sigma_1^{(\lambda^1)}]\otimes [\Sigma_2^{(\lambda^2)}]
   \big ).
\end{eqnarray*}
It follows that $\Psi \circ \big (\mathfrak p_{-i}^{X_1} \otimes
\text{\rm Id} \big ) = \mathfrak p_{-i}([\Sigma_1]) \circ \Psi$
for $i>0$.
\end{proof}

Let $\Lambda_1$ be the ring $\Lambda$ of symmetric functions with
the bilinear form
\begin{eqnarray*}
\langle p_{\lambda}, p_{\mu}\rangle = \delta_{\lambda,\mu}
\frac{1}{\mathfrak z_{\lambda}},
\end{eqnarray*}
$\Lambda_2$ be the same ring $\Lambda$ equipped with a different
bilinear form
\begin{eqnarray*}
\langle p_{\lambda}, p_{\mu}\rangle = \delta_{\lambda,\mu} \,\,
\frac{1}{\mathfrak z_{\lambda}} \,\, \big ({1}/({\gamma -1}) \big
)^{\ell(\lambda)},
\end{eqnarray*}
and $\Lambda_i^{n_i}$ be the space of degree-${n_i}$ symmetric
polynomials in $\Lambda_i$. The tensor product $\Lambda_1
\otimes_{\mathbb C}\Lambda_2$ has an induced bilinear form. Let
\begin{eqnarray*}
(\Lambda_1 \otimes_{\mathbb C}\Lambda_2)^n=\bigoplus_{n_1+n_2=n}
\Lambda^{n_1}_1\otimes_{\mathbb C}\Lambda_2^{n_2}.
\end{eqnarray*}
We define a ring structure on $(\Lambda_1\otimes_{\mathbb
C}\Lambda_2)^n$ by declaring the elements
\begin{eqnarray*}
\frac{P_{\lambda^1}^{(1)}}{c^{\prime}_{\lambda^1}(-1, -1)} \otimes
\frac{P_{\lambda^2}^{(1/(\gamma-1))}}
{c^{\prime}_{\lambda^2}(\gamma -1, 1)} \in
(\Lambda_1\otimes_{\mathbb C}\Lambda_2)^n,
\end{eqnarray*}
where $|\lambda^1| + |\lambda^2| =n$, to be idempotents.

\begin{theorem}\label{thm2}
There exists a linear isomorphism $\Phi\colon \mathbb H_X\to
\Lambda_1\otimes_{\mathbb C}\Lambda_2$ preserving bilinear forms
such that $\Phi \big ( \mathfrak p_{-\lambda^1}([\Sigma_1])
\mathfrak p_{-\lambda^2}([\Sigma_2])|0\rangle \big )
=p_{\lambda^1} \otimes p_{\lambda^2}$,
\begin{eqnarray*}
\Phi \big ([\Sigma_1^{(\lambda^1} \times \Sigma_2^{(\lambda^2)}]
\big )=m_{\lambda^1}\otimes m_{\lambda^2},
\end{eqnarray*}
and $\Phi([\lambda^1, \lambda^2])=P^{(1)}_{\lambda^1} \otimes
P^{(1/(\gamma-1))}_{\lambda^2}$. Furthermore, the restriction of
the map $\Phi$ to $\mathbb H_n$ is a ring isomorphism onto
$(\Lambda_1 \otimes_{\mathbb C}\Lambda_2)^n$.
\end{theorem}
\begin{proof}
Follows from Lemma~\ref{lem2} and arguments similar for
Theorem~\ref{thm1}.
\end{proof}

Next, we discuss the implication of Theorem~\ref{thm2} to the
ordinary cohomology ring $H^*(\Xn)$. With the notations for
Heisenberg algebras in the ordinary cohomology setting in
\cite{Na2}, a linear basis of $H^*(\Xn)$ consists of classes of
the form
\begin{eqnarray} \label{bas1}
\mathfrak Q_{\lambda^1,\lambda^2} \,\,\, \stackrel{\text{def}}{=}
\,\,\, \prod_{k \ge 1} P_{1_X}[-k]^{m_k(1)} \prod_{k \ge 1}
P_{\Sigma_0}[-k]^{m_k(2)} \vac
\end{eqnarray}
where $|\lambda^1| + |\lambda^2| =n$, $\lambda^i =(1^{m_1(i)}
2^{m_2(i)} \cdots)$, $1_X \in H^0(X)$ is the fundamental
cohomology class, and $\Sigma_0$ is the zero-section in $X =
X(\gamma)$. Denote
\begin{eqnarray} \label{Tbas1}
\mathfrak Q^{\mathbb T}_{\lambda^1,\lambda^2} \,\,\,
\stackrel{\text{def}}{=} \,\,\,  \prod_{k \ge 1} \mathfrak
p_{-k}(t)^{m_k(1)} \prod_{k \ge 1} \mathfrak
p_{-k}([{\Sigma_0}])^{m_k(2)} \vac \in \mathbb H_n.
\end{eqnarray}
The graded element in $H^*(\Xn)$ associated to $\mathfrak
Q^{\mathbb T}_{\lambda^1,\lambda^2}$ is $\mathfrak
Q_{\lambda^1,\lambda^2}$. Since $\{t, [{\Sigma_0}]\}$ and
$\{[\Sigma_1], [\Sigma_2]\}$ are two different linear basis for
$H_{\mathbb T}^2(X)$, it follows from the definitions that there
is a simple explicit transition matrix $M_1$ between the basis of
the equivariant cohomology $H^*_{\mathbb T} (\Xn)$ given by
(\ref{Tbas1}) and the basis given by
\begin{eqnarray} \label{Tbas2}
\mathfrak R^{\mathbb T}_{\lambda^1,\lambda^2} \,\,\,
\stackrel{\text{def}}{=} \,\,\, \prod_{k \ge 1} \mathfrak
p_{-k}([\Sigma_1])^{m_k(1)} \prod_{k \ge 1} \mathfrak
p_{-k}([\Sigma_2])^{m_k(2)} \vac.
\end{eqnarray}
For a fixed $r$, denote by $\big (g_{\lambda, \mu}^{(r)} \big )$
the transition matrix between Jack polynomials $P^{(r)}_\lambda$
and the power-sums $p_\mu$. By Theorem~\ref{thm2}, the transition
matrix between the basis (\ref{Tbas2}) and the basis
$\{[\mu^1,\mu^2]\}$ is provided by $M_2 \,\,
\stackrel{\text{def}}{=} \,\, \big (g_{\lambda^1, \mu^1}^{(1)}
\big ) \otimes \big (g_{\lambda^2, \mu^2}^{(1/(\gamma-1))} \big
)$.

Now the structure of the ordinary cohomology ring $H^*(\Xn)$ can
be described as follows in terms of the Heisenberg monomials
(\ref{bas1}). Given two Heisenberg monomials $\mathfrak
Q_{\lambda^1, \lambda^2}$ and $\mathfrak Q_{\mu^1, \mu^2}$ in
$H^*(\Xn)$, we consider the product of $\mathfrak Q^{\mathbb
T}_{\lambda^1, \lambda^2}$ and $\mathfrak Q^{\mathbb T}_{\mu^1,
\mu^2}$ in $\mathbb H_n$. The latter can be calculated by
transferring over to the basis $\{ \mathfrak R^{\mathbb
T}_{\lambda^1,\lambda^2}\}$ via the transition matrix $M_1$, and
then to the basis $\{[\mu^1,\mu^2]\}$ via the transition matrix
$M_2$, where the product is explicitly known. Finally, we reverse
the above steps and pass to the associated graded ring to express
the product of $\mathfrak Q_{\lambda^1, \lambda^2}$ and $\mathfrak
Q_{\mu^1, \mu^2}$ in terms of the Heisenberg monomial basis
(\ref{bas1}) of $H^*(\Xn)$.

\end{document}